\documentclass[12pt,a4paper,oneside]{amsart}

\usepackage{amssymb,amscd}
\usepackage[latin1]{inputenc}

\def\w{\omega}

\def\cal{\mathcal}

\def\proof{\noindent{\sl Proof}. }
\def\endproof{\hfill$\square$\smallskip\\}

\newtheorem{thm}{Theorem}
\newtheorem{cor}{Corollary}
\newtheorem{lem}{Lemma}
\newtheorem{prop}{Proposition}
\newtheorem{defn}{Definition}

\makeindex \makeglossary
\title{Generalized Vandermonde's system and Lagrange's interpolation}
\begin{document}
\maketitle
\begin{center}
{\sc Jean-Philippe PR\' EAUX}\footnote[1]{\noindent Research center of the French air force (CReA), F-13661 Salon de
Provence air, France}
\footnote[2]{Laboratoire d'Analyse Topologie et Probabilités,  Universit\'e de Provence, 39 rue F.Joliot-Curie, F-13453
marseille cedex 13, France.
\\ \indent\ {\it E-mail }: \ preaux@cmi.univ-mrs.fr}
, {\sc Jacques RAOUT}\footnote[3]{\noindent Morpho-Analysis in Signal processing lab., Research center of the French
air force (CReA), F-13661 Salon de Provence air, France.
\\ \indent {\it E-mail }: \ jraout@cr-ea.net\smallskip }

\end{center}

\begin{abstract} We give explicit formulas as well as a quadratic
time algorithm to solve (so called {\sl generalized Vandermonde's}) systems of linear equations with $p$ equations and
$n$ variables. It allows in particular to find all (so called {\sl Lagrange's interpolation}) polynoms with degree
$n-1$ taking given values in $p$ distinct given points.
\end{abstract}

\section*{Introduction}
Vandermonde's linear systems of equations (Alexandre Vandermonde, french Mathematician, 1735--1796) naturally appear in
numerical analysis  to find Lagrange's interpolation polynom. When one wants to determine a polynom $P(x)$ with degree
$n-1$ taking given values $q_1,q_2,\ldots ,q_n$ in $n$ distinct points $a_1,a_2,\ldots ,a_n$, one has to solve the (so
called {\sl Vandermonde's}) linear system $(*)$ of equations with unknowns $w_1,w_2,\ldots ,w_n$ :
$$\left\lbrace
    \begin{array}{l}
    \w_1+\w_2a_1+\w_3a_1^2+\cdots + \w_na_1^{n-1}=q_1\\
    \w_1+\w_2a_2+\w_3a_2^2+\cdots + \w_na_2^{n-1}=q_2\\
    \vdots\qquad \qquad\qquad\qquad\qquad\qquad\qquad\quad \vdots\\
    \w_1+\w_2a_n+\w_3a_n^2+\cdots + \w_na_n^{n-1}=q_n\\
    \end{array}\right. $$
 It admits a unique solution and the (so called {\sl Lagrange's interpolation})
polynom $P(x)$ of couples $\lbrace (a_1,q_1),$ $(a_2,q_2),\ldots ,(a_n,q_n)\rbrace$ is given by the formula :
$$P(x)=\sum_{k=1}^n\w_k x^{k-1}$$

Vandermonde's linear systems also appear naturally in several problems of linear algebra.\smallskip\\

The aim of this work is to give explicit solution to the more general interpolation problem consisting in finding all
polynoms of degree $n-1$ with given values in $p$ distinct points, where $p$ may be different from $n$. It brings us to
any system of $p$ equations as above with $n$ unknowns, that we shall extensively call {\sl generalized Vandermonde's
linear systems}.

On the one hand we provide in \S 1 explicit formulas solving the problem in full generality, and on the other we build
in \S 2 a quadratic time algorithm implementing the solution.

\section{Explicit solutions to generalized Vandermonde's linear system} We first consider in \S 1.1 preliminaries notions necessary
for effective computations. In \S 1.2 we establish explicit formulas for the inverse of a Vandermonde's square matrix,
and consequently for the Lagrange's interpolation polynom. The two remaining sections are concerned with generalized
Vandermonde's system : in \S 1.3 we explicit the kernel of a non square Vandermonde's matrix, and in \S 1.4 we collect
all preceding result to give a solution to the general problem.

\subsection{Preliminaries : monomial coefficients} In this section, $A$ stands for a unitary commutative ring.
Given $p$ distinct elements $a_1,a_2,\ldots ,a_p$ in $A$, we define for any positive integer $t$ the {\sl monomial
coefficient of codegree} $t$ on $a_1,a_2,\ldots ,a_p$ : $\sigma_{a_1\cdots a_p}(t)$ (or more concisely $\sigma(t)$ when
there is no ambiguity) informally speaking as the sum of all products of $t$-uples on $a_1,a_2,\ldots ,a_p$ without
repetition.

\begin{defn} Let $a_1,a_2,\ldots ,a_p$, be $p$ distinct elements in $A$. For any positive integer $t$, the
{\sl monomial coefficient of codegree} $t$ on $a_1,a_2,\ldots ,a_p$ is the element  $\sigma_{a_1\cdots a_p}(t)$ of $A$
(or more concisely $\sigma(t)$) defined by :
$$
\left\lbrace\begin{array}{lll} \sigma(0)
&=1_A &\\
 \sigma(t)&=\underset{a_1\cdots
a_p}{\sigma}(t)=\sum_{1\le i_1< \cdots< i_t\le p}
a_{i_1}a_{i_2}\cdots a_{i_t}\quad &\text{$\forall\, t=1,2,\ldots ,p$}\\
\sigma(t)&=\underset{a_1\cdots a_p}{\sigma}(t)=0_A &\text{$\forall\, t>p$}$$
\end{array}\right.
$$
\end{defn}

\indent For more convenience in transcription of formulas, we also define the related notation :

\begin{defn} Under the same hypothesis as above, and given $a=a_k$ with $1\le k\le p$,
the $\widehat{a}$-{\sl monomial coefficient of codegree} $t$ on $a_1,a_2,\ldots ,a_p$, is defined as the {\sl monomial
coefficient of codegree} $t$ on $a_1,\ldots ,\widehat{a_k},\ldots ,a_p$, that is :
$$
\overline{\sigma}^a(t)=\underset{a_1\cdots a_p}{\overline{\sigma}^a}(t)=\sum_{1\le i_1< \cdots< i_t\le p\atop
i_1,\ldots , i_t\not=k} a_{i_1}a_{i_2}\cdots a_{i_t}
$$
\end{defn}

We now explicit basic properties of monomial coefficients.
\begin{lem}\label{sigmainduction}
Under the same hypothesis as above, for any $t=0,1,\ldots, p$, one has :
$$\sigma(t)=\overline{\sigma}^a(t)+a\,\overline{\sigma}^a(t-1)$$
\end{lem}
\noindent {\sl Proof.} Immediate by definition \endproof

We now enonce the fundamental property of monomial coefficients which motivates them to be introduced :

\begin{prop}\label{equation}
In the polynom ring $A[x]$, one verifies :
$$\prod_{i=1}^p (x-a_i)=\sum_{i=0}^{p} (-1)^{p-i}
\underset{a_1\cdots a_p}{\sigma}(p-i)\, x^i=\sum_{i=0}^{\infty} (-1)^{p-i} \underset{a_1\cdots a_p}{\sigma}(p-i)\,
x^i$$
\end{prop}
\proof Immediate by induction
\endproof

\noindent One immediately obtains by setting $x=a$ :
\begin{cor}
Under the same hypothesis as above one has for any $a=a_1,\ldots ,a_p$ :
$$\sum_{i=0}^{p}(-1)^ia^i\sigma(p-i)=0$$
\end{cor}

\subsection{Inversion of a square Vandermonde's matrix}
Consider the linear system $(*)$ as appearing in the introduction.
Let $a_1,a_2,\ldots , a_n$ be $n$ elements of a commutative field. The {\sl Vandermonde's square matrix} is the matrix
$n\times n$ defined by :
$$V_n=\left(\begin{array}{ccccc}
1 & a_1 & a_1^2 & \cdots & a_1^{n-1}\\
1 & a_2 & a_2^2 & \cdots & a_2^{n-1}\\
\vdots & & & &\vdots \\
1 & a_n & a_n^2 & \cdots & a_n^{n-1}\\
\end{array}\right)$$
 \noindent One easily verifies that its determinant is given by :
$$\det(V_n)=\prod_{1\leq j<i\leq n}^n (a_i-a_j)\not=0$$
so that the matrix $V_n$ is inversible.\smallskip\\
\indent Consider the family of $n$ polynoms defined for all $j=1,2,\ldots ,n$ by :
$$P_j(x)=\prod_{k=1\atop k\not= j}^n\frac{x-a_k}{a_j-a_k}=\sum_{i=1}^n c_{i,j}\, x^{i-1}$$
Clearly :
$$
P_j(a_i)=\left\lbrace\begin{array}{l} 1\quad\text{if $i=j$}\\
0\quad\text{otherwise}
\end{array}\right.
$$
so that the inverse matrix $V_n^{-1}$ of $V_n$ is given by :
$$V_n^{-1}=\left(c_{i,j}\right)_{1\leq i\leq n\atop 1\leq j \leq n}=
\left(\begin{array}{cccc}
c_{1,1} & c_{1,2} & \cdots & c_{1,n}\\
c_{2,1} & c_{2,2} & \cdots & c_{2,n}\\
\vdots & &  &\vdots \\
c_{n,1} & c_{n,2} & \cdots & c_{n,n}\\
\end{array}\right)
$$
with the $c_{i,j}$ defined as above.

Let us denote by $D_j$ the denominator of $P_j(x)$, it lies in $K^*$, and by $N_j(x)$ the numerator of $P_j(x)$, it
lies in $K[x]$ and has degree $n-1$. It follows from proposition \ref{equation} that :
$$N_j(x)= \prod_{k=1\atop k\not= j}^n  (x-a_k) =
\sum_{i=1}^{n} (-1)^{i-1}\, \overline{\sigma}^{a_j}(n-i)\, x_j^{i-1}$$
$$D_j=\prod_{k=1\atop k\not= j}^n (a_j-a_k)=\sum_{i=1}^{n} (-1)^{i-1}\, \overline{\sigma}^{a_j}(n-i)\,
a_j^{i-1}$$ et donc pour tout $i,j=1,2,\ldots ,n$
$$c_{i,j}  =\frac{ \overline{\sigma}^{a_j}(n-i)}
{\sum_{k=1}^{n} (-1)^{i+k}\ \overline{\sigma}^{a_j}(n-k)\, a_j^{k-1}}$$

Hence, the system $(*)$  :
$$
\left(\begin{array}{ccccc}
1 & a_1 & a_1^2 & \cdots & a_1^{n-1}\\
1 & a_2 & a_2^2 & \cdots & a_2^{n-1}\\
\vdots & & & &\vdots \\
1 & a_n & a_n^2 & \cdots & a_n^{n-1}\\
\end{array}\right)
\left(\begin{array}{c}
\w_1\\
\w_2\\
\vdots\\
\w_n
\end{array}\right)=
\left(\begin{array}{c}
q_1\\
q_2\\
\vdots\\
q_n\end{array}\right)
$$
has a unique $n$-uple solution $(\w_1,\w_2,\ldots ,\w_n)$, given by, $\forall\, i=1,2,\ldots ,n$ :
$$\w_i=\sum_{j=1}^n\ \frac{ \overline{\sigma}^{a_j}(n-i)\, q_j}
{\sum_{k=1}^{n} (-1)^{i+k}\ \overline{\sigma}^{a_j}(n-k)\, a_j^{k-1}}
$$

One can explicit the Lagrange's interpolation polynom with degree $n-1$ and respective values $q_1,q_2,\ldots ,q_n$ in
points $a_1,a_2,\ldots ,a_n$ :
$$P(x)=\sum_{i=1}^n \w_i\, x^{i-1}=\sum_{i=1}^n\ \sum_{j=1}^n\ \frac{ \overline{\sigma}_{a_1\cdots a_n}^{a_j}(n-i)\, q_j}
{\sum_{k=1}^{n} (-1)^{i+k}\ \overline{\sigma}_{a_1\cdots a_n}^{a_j}(n-k)\, a_j^{k-1}}\ x^{i-1}
$$

\subsection{Kernel of a Vandermonde's matrix $p\times n$} In this section $p$ and $n$ are integers with $1\le p\le n$ and
$a_1,a_2,\ldots ,a_p$ are distinct elements in a commutative field $K$. One considers the {\sl Vandermonde's matrix
with $p$ lines and $n$ columns} :
$$V_{p,n}={\left(\begin{array}{ccccc}
1 & a_1 & a_1^2 & \cdots & a_1^{n-1}\\
1 & a_2 & a_2^2 & \cdots & a_2^{n-1}\\
\vdots & & & & \\
1 & a_p & a_p^2 & \cdots & a_p^{n-1}\\
\end{array}\right)}_{p\times n}$$

The following theorem explicits a basis for the kernel of $V_{p,n}$.
\begin{thm}
The kernel of $V_{p,n}$ has dimension $n-p$ and admits a basis given by the vector $\overset{\rightarrow}{v}_1$
together with all its cyclic conjugates $\overset{\rightarrow}{v}_2,\ldots, \overset{\rightarrow}{v}_{n-p}$ :
$$\overset{\rightarrow}{v}_1=\left(\begin{array}{c}
(-1)^p\sigma(p)\\
\vdots\\
(-1)^{i}\sigma(i)\\
\vdots\\
-\sigma(1)\\
 1\\
 0\\
 \vdots\\
 0\\
 0\\
 \end{array}\right)\ ;\
\overset{\rightarrow}{v}_2=\left(\begin{array}{c}
0\\
(-1)^p\sigma(p)\\
\vdots\\
(-1)^{i}\sigma(i)\\
\vdots\\
-\sigma(1)\\
 1\\
 0\\
 \vdots\\
 0\\
 \end{array}\right)\ ;\cdots ;\
\overset{\rightarrow}{v}_{n-p}=\left(\begin{array}{c}
0\\
\vdots\\
\vdots\\
0\\
(-1)^p\sigma(p)\\
\vdots\\
(-1)^{i}\sigma(i)\\
\vdots\\
-\sigma(1)\\
 1\\
  \end{array}\right)
 $$
\end{thm}

\proof One can easily establish that $V_{p,n}$ has maximal rank. Hence, $\ker V_{p,n}$ has dimension $n-p$. It follows
from corollary 1 that $v_1$ lies in $\ker V_{p,n}$. By multiplying equation of corollary 1 by $a^k$, for $0<k<n-p$, one
sees that $v_2,v_3,\ldots,v_{n-p}$ also lie in $\ker V_{p,n}$. Moreover the family of vectors $v_1,v_2,\ldots, v_{n-p}$
is clearly free, so that they constitue a basis of $\ker V_{p,n}$.\endproof \indent (Note that the kernel of the
Vandermonde's matrix with $n$ lines and $p$ columns $V_{n,p}$ has obviously dimension 0.)

\subsection{Generalized Vandermonde's system} Collecting all the above :

\begin{thm}
The generalized Vandermonde's linear system with $p$ equations and $n$ unknowns :
$$
\left(\begin{array}{ccccc}
1 & a_1 & a_1^2 & \cdots & a_1^{n-1}\\
1 & a_2 & a_2^2 & \cdots & a_2^{n-1}\\
\vdots & & & & \\
1 & a_p & a_p^2 & \cdots & a_p^{n-1}\\
\end{array}\right)_{p,n}
\left(\begin{array}{c}
\w_1\\
\w_2\\
\vdots\\
\w_n
\end{array}\right)=
\left(\begin{array}{c}
q_1\\
q_2\\
\vdots\\
q_p\end{array}\right)
$$
has solutions space $\cal S$ a codimension $p$ affine sub-space of the $n$-dimensional $K$-vector space $K^n$ : $\cal
S=\overset{\rightarrow}{\w}_0+\ker(V_{p,n})$, with $\overset{\rightarrow}{\w}_0^t=(\w_1,\w_2,\ldots, \w_p,0,\ldots,0)$,
and
$$
                            \left(\begin{array}{c}
                                    \w_1\\
                                    \w_2\\
                                    \vdots\\
                                    \w_p\\
                                    \end{array}\right)
                                    =
                                    \left(\begin{array}{ccccc}
1 & a_1 & a_1^2 & \cdots & a_1^{p-1}\\
1 & a_2 & a_2^2 & \cdots & a_2^{p-1}\\
\vdots & & & & \\
1 & a_p & a_p^2 & \cdots & a_p^{p-1}\\
\end{array}\right)_{p,p}
                                     \left(\begin{array}{c}
                                    q_1\\
                                    q_2\\
                                    \vdots\\
                                    q_p\\
                                    \end{array}\right)
$$
\end{thm}


\subsection{Example}
Let :
$$V_{n-1,n}={\left(\begin{array}{ccccc}
1 & a_2 & a_2^2 & \cdots & a_2^{n-1}\\
1 & a_3 & a_3^2 & \cdots & a_3^{n-1}\\
\vdots & & & & \\
1 & a_n & a_n^2 & \cdots & a_n^{n-1}\\
\end{array}\right)}$$

$$\ker
V_{n-1,n}=Vect(
    \left(
        \begin{array}{l}
        (-1)^{n-1}\Sigma(n-1)\\
        \vdots\\
        (-1)^i\Sigma(i)\\
        \vdots\\
        -\Sigma(1)\\
        1
        \end{array}\right)
        )
        =Vect(
        \left(
        \begin{array}{l}
        (-1)^{n-1}a_2a_3\cdots a_n\\
        \vdots\\
        (-1)^i\Sigma(i)\\
        \vdots\\
        -(a_2+a_3+\cdots+a_n)\\
        1
        \end{array}\right)
        )$$

\section{Algorithmic computations}
All along the section a non negative integer $p$ as well as $p$ elements $a_1,a_2,\ldots , a_p$ in a unitary
commutative ring $A$ are given.

\subsection{Computation of the monomial coefficients}
 The lemma
\ref{sigmainduction} provides an algorithm {\it à la} Pascal  to compute all monomial coefficients $\sigma(1),
\sigma(2)$, $\ldots ,\sigma(p)$ :

\begin{prop} There is an algorithm with $O(p^2)$ complexity which computes all monomial coefficients ${\sigma}(t)$,
for $t=0,1,\ldots ,p$ :

Consider a matrix with $p+1$ rows and $p$ columns, where initially the left column is constitued of $1$ and on the top
row is followed with zeros. Then fill in the remaining elements starting from the top row to the bottom one and from
left to  right using :
$$\text{element$(i,j)$}= \text{element$(i-1,j)+a_i.$element$(i-1,j-1)$}$$
At the end of the process the bottom row consists in $\sigma(0),\sigma(1),\ldots , \sigma(p)$.
\end{prop}

\noindent{\bf Example.} Given three elements $a,b,c$, we apply this method :\\

\centerline{
\begin{tabular}{r||l|l|l|l}
  & 0 & 1 & 2 & 3\\
\hline \hline & 1 & 0 & 0 & 0\\
\hline $\sigma_a$ & 1 & a & 0 &0\\
\hline $\sigma_{ab}$ & 1 & a+b & ab &0\\
\hline $\sigma_{abc}$ & 1 & a+b+c & ab+ac+bc & abc\\
\end{tabular}
}
\medskip
 \noindent Successive computations give :\smallskip\\ ${\sigma}_a(1)=0+a\times 1=a$\smallskip\\
${\sigma}_{ab}(1)={\sigma}_a(1)+b\,{\sigma}_a(0)=a+b\times 1=a+b$ \\
 ${\sigma}_{ab}(2)={\sigma}_a(2)+b\,{\sigma}_a(1)=0+b\times a=ab$\smallskip\\
${\sigma}_{abc}(1)={\sigma}_{ab}(1)+c\,{\sigma}_{ab}(0)=a+b+c\times 1=a+b+c$ \\
${\sigma}_{abc}(2)={\sigma}_{ab}(2)+c\,{\sigma}_{ab}(1)=ab+c\times (a+b)=ab+ac+bc$ \\
${\sigma}_{abc}(3)={\sigma}_{ab}(3)+c\,{\sigma}_{ab}(2)=0+c\times ab=abc$\\

\noindent {\bf Algorithm.} Given the integer $p\geq 1$ and  $p$ distinct numbers $a_1,a_2,\ldots ,a_p$, the following
algorithm, written in classical algorithmic language, computes all $\sigma(t)$ for $t=0,1,\ldots ,p$ :\medskip\\
{\tt const\ integer\qquad\qquad\quad p\qquad \qquad\qquad\qquad\qquad { \verb+//+ \sl contains $p$}\\
const\ array of number\quad A[p]=[a1,a2,$\ldots$ , ap]\ { \verb+//+ \sl contains $a_1,a_2,\ldots ,a_p$}\\
array of number\qquad\qquad S[p+1]=[1,0,$\ldots$ ,0]\\
{\verb+//+ \sl First element of an array has index 0 !}\\
for i=1 to p\\
\indent\indent for j=i to 1 step -1\\
\indent\indent\indent\indent S[j]=S[j]+A[i-1]*S[j-1]\\
{\verb+//+ \sl {\tt S[p+1]} contains $\Sigma(0),\Sigma(1),\ldots,
\Sigma(p)$.}\bigskip\\
 }

Lemma 1 also provides an algorithm to compute all coefficients ${\overline{\Sigma}}^{a}(t)$ for all  $a=a_1,a_2,\ldots
,a_p$ and $t=0,1,\ldots ,p-1$.\\

\begin{prop} Once all monomial coefficients ${\sigma}(0), {\sigma}(1),\ldots
,{\sigma}(p)$, as well as $a=a_i$ are given, one can algorithmically compute the coefficients $\overline{\sigma}^a(0),
\overline{\sigma}^a(1), \ldots ,\overline{\sigma}^a(p-1)$ in linear time $O(p)$, using the inductive formula :
$$\left\lbrace\begin{array}{l}
                \overline{\sigma}^a(0)=1\\
                \overline{\sigma}^a(n)=\sigma(n)-a\,\overline{\sigma}^a(n-1)
                \end{array}\right.
                $$
En particular, using proposition 2, given $p$ distinct elements $a_1,a_2,\ldots ,a_p$, there is an algorithm which
returns in quadratic time $O(p^2)$ the sequence $\sigma(1),\sigma(2),\ldots ,\sigma(p)$ as well as  $p$ sequences
$\overline{\sigma}^{a_i}(1), \overline{\sigma}^{a_i}(2),\ldots ,\overline{\sigma}^{a_i}(p-1)$, for $i=1, 2,\ldots, p$.
\end{prop}

\noindent{\bf Example.} Given three distinct elements  $a,b,c$, we apply this method :\smallskip\\

\centerline{
\begin{tabular}{r||l|l|l|l}
 & 0 & 1 & 2 & 3\\
\hline
\hline $\sigma^{\phantom{a}}$ & 1 & a+b+c & ab+ac+bc & abc\\
\hline $\overline{\sigma}^a$ & 1 & b+c & bc &0\\
\hline $\overline{\sigma}^b$ & 1 & a+c & ac &0\\
\hline $\overline{\sigma}^c$ & 1 & a+b & ab &0\\
\end{tabular}
}

\noindent Soit :\smallskip\\
$\overline{\sigma}^a(1)={\sigma}(1)-a\,\overline{\sigma}^a(0)=a+b+c - a\times 1=b+c$ \\
$\overline{\sigma}^a(2)={\sigma}(2)-a\,\overline{\sigma}^a(1)=ab+ac+bc-a(b+c)=bc$\smallskip\\
$\overline{\sigma}^b(1)={\sigma}(1)-b\,\overline{\sigma}^b(0)=a+b+c - b\times 1=a+c$ \\
$\overline{\sigma}^b(2)={\sigma}(2)-b\,\overline{\sigma}^b(1)=ab+ac+bc-b(a+c)=ac$\smallskip\\
$\overline{\sigma}^c(1)={\sigma}(1)-c\,\overline{\sigma}^c(0)=a+b+c - c\times 1=a+b$ \\
 $\overline{\sigma}^c(2)={\sigma}(2)-c\,\overline{\sigma}^c(1)=ab+ac+bc-c(a+b)=ab$\smallskip\\

\noindent{\bf Algorithm.} Given the integer $p\geq 1$,  $p$ numbers  $a_1,a_2,\ldots ,a_p$ and the sequence of monomial
coefficients $\sigma(0),\sigma(1),\ldots ,\sigma(p)$, the following algorithm written in standard algorithmic language
returns all $\widehat{a}$-monomial coefficients $\overline{\sigma}^a(t)$
for $t=0,1,\ldots ,p-1$ and $a=a_1,a_2,\ldots ,a_p$. \\

{\tt
\noindent const\quad integer\qquad\qquad\quad\ p\qquad\qquad\qquad\qquad\qquad{ \verb+//+ \sl contains $p$}\\
const\quad array of number\quad\ A[p]=[a1,\ldots ,ap] \qquad{ \verb+//+ \sl contains $a_1,a_2,\ldots ,a_p$}\\
const\quad array of number\quad\ S[p+1]=[s0,s1,\ldots ,sp]\ {\verb+//+ $\sigma(0),\sigma(1),\ldots ,\sigma(p)$}\\
array of array of number\ T[p][p]\smallskip\\
for i=0 to p-1\\
\indent\indent T[i][0]=1\\
\indent\indent for j=1 to p-1\\
\indent\indent\indent\indent T[i][j]= S[j]-a*T[i][j-1]\smallskip\\
\verb+\\+  {\tt T[i]} {\sl contains\ $\overline{\sigma}^{a_i}(0), \overline{\sigma}^{a_i}(1), \ldots ,
\overline{\sigma}^{a_i}(p-1)$ }}

\end{document}